\newif\ifsmfart
\IfFileExists{smfart.cls}
   {\documentclass[12pt,english]{smfart}
    \IfFileExists{smfenum.sty}{\usepackage{smfenum}}{}
    \usepackage{bull}
    \smfarttrue}
{\message{^^J*** It would be better to typeset this file with smfart.cls ***^^J^^J}

\documentclass[12pt]{amsart}}
\usepackage{amssymb}
\numberwithin{equation}{section}

\ifx\pdftexversion\undefined
  \usepackage[dvips]{graphicx}
\else
  \usepackage[pdftex]{graphicx}
\fi

\input xy
\xyoption{all}

\advance\headheight 2pt
\calclayout
\allowdisplaybreaks[3]

\theoremstyle{plain}
\newtheorem{prop}{Proposition}

\newtheorem{theo}[prop]{Theorem}

\newtheorem{lemm}[prop]{Lemma}

\theoremstyle{definition}
\newtheorem{defi}[prop]{Definition}

\newtheorem{conj}[prop]{Conjecture}

\theoremstyle{remark}

\newtheorem{exam}[prop]{Example}

\newcommand{\bA}{\mathbb A}
\newcommand{\bC}{\mathbb C}
\newcommand{\bN}{\mathbb N}
\newcommand{\bP}{\mathbb P}

\newcommand{\cL}{\mathcal L}
\newcommand{\cN}{\mathcal N}
\newcommand{\cO}{\mathcal O}

\newcommand{\cU}{\mathcal U}
\newcommand{\cV}{\mathcal V}

\newcommand{\cX}{\mathcal X}
\newcommand{\cY}{\mathcal Y}
\newcommand{\ra}{\rightarrow}
\newcommand{\lra}{\longrightarrow}
\newcommand{\Spec}{{\rm Spec}}
\newcommand{\Hom}{{\rm Hom}}
\newcommand{\iHom}{{\mathcal H}\!{\mathit om}}
\newcommand{\Ext}{{\rm Ext}}
\newcommand{\iExt}{{\mathcal E}\!{\mathit xt}}
\newcommand{\mfm}{\mathfrak m}
\newcommand{\mfp}{\mathfrak p}
\newcommand{\mfq}{\mathfrak q}
\newcommand{\mfr}{\mathfrak r}

\makeatother
\makeatletter

\def\Title    {Weak approximation over function fields}
\def\Author   {Brendan Hassett and Yuri Tschinkel}
\def\Subject  {Algebraic geometry}
\def\Keywords {rationally connected varieties, weak approximation}

\newif\ifpdf
   \ifx\pdfoutput\undefined
   \pdffalse
   \else
           \pdfoutput=1
            \pdftrue
   \fi
      \ifpdf
\usepackage[pdftex,
hyperindex=true,
pdftitle={\Title},pdfauthor={\Author},pdfsubject={\Subject},
pdfkeywords={\Keywords},pdfpagemode={UseOutlines},
bookmarks=true,bookmarksopen=true,
bookmarksopenlevel=0,bookmarksnumbered=true,
pdfstartview={FitV},
colorlinks=true
]{hyperref}
 \else
 \usepackage{hyperref}
 \fi

\begin{document}

\title[Weak approximation over function fields]
{Weak approximation over function fields}
\author{Brendan Hassett}
\address{Department of Mathematics \\
Rice University, MS 136 \\
Houston, TX 77251-1892}
\email{hassett@rice.edu}
\author{Yuri Tschinkel}
\address{Mathematisches Institut\\
Bunsenstr. 3-5 \\
37073 G\"ottingen, Germany}
\email{yuri@uni-math.gwdg.de}

\date{\today}

\begin{abstract}
We prove that
rationally connected varieties over the function field of a complex curve
satisfy weak approximation for places of good reduction.
\end{abstract}

\maketitle
\tableofcontents

\section{Introduction}
\label{sect:introduction}

Let $F$ be a number field and $X$ an algebraic variety over $F$.
Does there exist an $F$-rational point on $X$?  If so, are they
ubiquitous on $X$? For many classes of varieties, such problems
are analyzed using local-to-global principles. 
The Hasse principle says that $X$ has an $F$-rational point
provided it has a rational point over each completion of $F$. 
The principle of weak approximation says that, given a finite collection of 
places of $F$ together with a point of $X$ over each of the corresponding 
completions, there exists an $F$-rational point approximating these
arbitrarily closely.

\

The impetus for this paper was 
the following result by Graber, Harris and Starr:  

\begin{theo}[\cite{GHS},Theorem 1.2]
\label{theo:great}
Let $F$ be the function field of a smooth curve over $\bC$. 
Every proper rationally connected variety $X$ over 
$F$ has an $F$-rational point. 
\end{theo}

An algebraic variety is rationally connected if any two points can be
joined by a rational curve (see Section~\ref{sect:rat} for more details).
Rational and unirational varieties are rationally connected. 
We refer the reader to \cite{dJS} for related results 
in positive characteristic.

\

From an arithmetic viewpoint, proving such a theorem 
entails surmounting two obstacles: First 
one needs to show that there are no obstructions to the existence 
of a local point, i.e., $X(F_{\nu})\neq \emptyset$ for all completions $F_{\nu}$
of $F$.  If $B$ is the smooth projective curve with $F=\bC(B)$
and $\pi:\cX\ra B$ is a model for $X$ over $B$,  
one has to show there exist local analytic sections of $\pi$
at each point of $B$.  
Secondly, one has to prove the Hasse principle for $X$ over $F$,
which entails constructing a global section of $\pi$. 
Note that over a number field, the Hasse principle may fail
even for cubic surfaces.

\

Theorem~\ref{theo:great} naturally leads one to ask 
whether rationally connected varieties over $\bC(B)$ satisfy weak 
approximation as well.  In this paper we prove this away from singular 
fibers of $\pi$, i.e., away from the places of bad reduction:
\begin{theo}
\label{theo:main}
Let $X$ be a smooth, proper,
rationally connected variety over the 
function field of a curve over $\bC$.  Then $X$ satisfies 
weak approximation for places of good reduction.
\end{theo}
Theorem~\ref{theo:great} and Theorem 2.13 of \cite{KMM} give
the zeroth-order case: 
There exists a section of $\pi:\cX \ra B$ passing through arbitrary points
of smooth fibers.

For varieties over function fields of curves, weak approximation 
is satisfied in the following cases \cite{CT}:
\begin{itemize}
\item stably rational varieties;
\item connected linear algebraic groups and homogeneous spaces
for these groups;
\item homogeneous space fibrations over varieties that satisfy weak approximation, 
for example, conic bundles over rational varieties;
\item Del Pezzo surfaces of degree at least four.
\end{itemize} 
Weak approximation is not known for general cubic surfaces. 
Madore has a manuscript addressing weak approximation for
cubic surfaces away from places of bad reduction.    

\

{\bf Acknowledgments:}
The first author was partially supported by the
Sloan Foundation and NSF Grants 0134259 and 0196187.
The second author was partially supported by NSF Grant 0100277.
Part of this work was done while both authors were visiting
the American Institute of Mathematics in Palo Alto. 
We thank J. de Jong, T. Graber, J. Harris, J. Koll\'ar and J. Starr
for conversations about this topic.

\section{Basic properties of weak approximation}
\label{sect:not}

\subsection{Definition}

Let $F$ be a number field or a function field
of a smooth projective curve $B$ over 
an algebraically closed ground field $k$ of
characteristic zero.  For each place $\nu$ of $F$,
let $F_{\nu}$ denote the $\nu$-adic completion of $F$.
Let $X$ be an algebraic variety 
of dimension $d$ over $F$;  in this paper,
all varieties are assumed to be geometrically integral.
Let $X(F)$ denote the set of $F$-rational points of $X$. 
One says that rational points on $X$ {\em satisfy weak approximation}
if, for any finite set of places 
$\{\nu_i\}_{i\in I}$ of $F$ and 
$\nu_i$-adic open subsets $U_i\subset X(F_{\nu_i})$, 
there is a rational point $x\in X(F)$
such that its image in each $X(F_{\nu_i})$ is contained in $U_i$.
In particular, for any collection of $x_i \in X(F_{\nu_i}), i\in I,$
there exists an $x\in X(F)$ arbitrarily close to each $x_i$.  

It is well known that weak approximation is a birational property:
If $X_1$ and $X_2$ are smooth varieties birational over $F$
then $X_1$ satisfies weak approximation if and only if $X_2$
satisfies weak approximation.  Given a smooth proper variety $X_1$, 
after applying
Chow's lemma and resolution of singularities we obtain a smooth
projective variety $X_2$ birational to $X_1$.  Thus in proving
weak approximation, it usually suffices to consider
projective varieties.  In particular, Theorem~\ref{theo:main}
reduces to this case.

\

For the rest of this paper we restrict our attention to the function
field case.  Places $\nu$ of $F$ correspond to points $b$ on $B$.
We also assume that $X$ is projective, so it admits a projective model
$$\pi:\cX \ra B,$$ 
i.e., a flat projective morphism 
with generic fiber $X$;
for each $b\in B$,
the fiber over $b$ is denoted
$$\cX_b=\cX\times_B \Spec(\cO_{B,b}/\mfm_{B,b}).$$
Sections of $\pi$ yield $F$-valued points of $X$ and 
conversely, each $F$-valued point of $X$ extends to a section of 
$\pi$.   
Let $\widehat{B}_b$ denote the completion of $B$ at $b$ and
$(\widehat{\cO}_{B,b},\widehat{\mfm}_{B,b})$ the 
associated complete local ring, which has fraction field $F_{\nu}$.
Restricting to this formal neighborhood of $b$ gives 

\centerline{
\xymatrix{
  \cX\times_B \widehat{B}_b \ar[d]_{\hat{\pi}_b} \ar[r] & \cX\ar[d]_{\pi} \\
  \widehat{B}_b \ar[r] & B 
 }}

\noindent
Sections of $\hat{\pi}_b$ restrict to $F_{\nu}$-valued points of $X$ and 
conversely, each $F_{\nu}$-valued point of $X$ extends to a section of 
$\hat{\pi}_b$.  Basic $\nu$-adic open subsets of $X(F_{\nu})$
consist of those sections of $\hat{\pi}_b$ 
which agree with a given section over 
$\Spec(\widehat{\cO}_{B,b}/\widehat{\mfm}^{N+1}_{B,b})
\subset\widehat{B}_b$, for
some $N\in \bN$. 
Weak approximation means that for any finite set of points $\{b_i\}_{i\in I}$
in $B$, sections $\hat{s}_i$ of $\hat{\pi}_{b_i}$,
and $N\in \bN$, there exists a section $s$ of $\pi$ 
agreeing with $\hat{s}_i$ modulo $\widehat{\mfm}^{N+1}_{B,b_i}$ for each $i$.

\subsection{Fibers of good reduction}
\label{sect:good}
We continue to assume that $X$ is a smooth projective variety over $F=k(B)$,
$\nu$ a place of $F$, and $b\in B$ the corresponding point.
A place $\nu$ is {\em of good reduction} for $X$ if there exists
a scheme
$$\widehat{\cX}_b \ra \widehat{B}_b,$$
proper and smooth over $\widehat{B}_b$,
with generic fiber isomorphic to $X$ over $F_{\nu}$.  
Let $S$ denote the finite set of places of bad reduction.

\begin{defi}
A variety 
$X$ satisfies {\em weak approximation for places of good reduction} if,
for any finite set of places of good reduction
$\{\nu_i\}_{i\in I}$ and
$\nu_i$-adic open subsets $U_i\subset X(F_{\nu_i})$,
there is a rational point $x\in X(F)$
such that its image in each $X(F_{\nu_i})$ is contained in $U_i$.
\end{defi} 
Suppose we have a model $\pi: \cX\ra B$ smooth over $B\setminus S$.  
Then we can express this in more geometric terms:
For each finite set of points $\{b_i\}_{i\in I}$
in $B\setminus S$, sections $\hat{s}_i$ of $\hat{\pi}_{b_i}$
and $N\in \bN$, there exists a section $s$ of $\pi$ 
agreeing with $\hat{s}_i$ modulo $\widehat{\mfm}^{N+1}_{B,b_i}$ for each $i$.

\begin{prop}\label{prop:goodmodel}
Retain the notation introduced above. 
There exists an algebraic space 
$$\pi:\cX \ra B,$$
proper and flat over $B$, smooth over the places of good reduction,
and with generic fiber $X$.  
Such a space is called a good model of $X$ over $B$.  
\end{prop}
\begin{proof}
Choose a projective model $\pi':\cX'\ra B$ for $X$ over $B$.
If $\pi'$ is smooth over $B\setminus S$ there is nothing to prove.  Otherwise,
let $\{b_j\}\subset B\setminus S$ denote the points over which 
$\cX'_{b_j}$ is singular; let $\widehat{\cX}'_{b_j}$ denote the 
completion of $\cX'$ along the central fiber $\cX'_{b_j}$.   
By assumption, 
there exists a proper smooth scheme
$$\widehat{\cX}_{b_j}\ra \widehat{B}_{b_j},$$
which is isomorphic to
our original model over the generic point.
Resolving the indeterminacy of the rational map
$$\widehat{\cX}_{b_j}\dashrightarrow \widehat{\cX}'_{b_j},$$
we find that these are
related by a sequence of modifications in the central fiber.  
This gives a sequence of formal modifications 
to $\cX'$ along the singular fibers, 
in the sense of Artin \cite{Art70}.
Theorems 3.1 and 3.2
of \cite{Art70} give a unique proper algebraic
space $\pi:\cX \ra B$ realizing these formal modifications to $\cX'$.  
\end{proof}
\begin{exam}
There are simple examples justifying the introduction of algebraic spaces.
Let $\pi:\cX'\ra B$ be a flat projective morphism such that each
fiber is a cubic surface with rational double points and the
generic fiber $X$ is smooth.  Suppose that near each point $b\in B$ the
local monodromy representation
$$
\mathrm{Gal}(\bar{F}_{\nu}/F_{\nu})\ra 
\mathrm{Aut}(\mathrm{Pic}(\bar{X}))
$$
is trivial.  
By a theorem of Brieskorn \cite{Br1} \cite{Br2}, 
there exists a simultaneous resolution
$$
\begin{array}{rcccl}
\cX & & \stackrel{\varrho}{\lra} & &\cX' \\
    & \searrow & &  \swarrow&  \\
    & 			 &       B & & 
\end{array},
$$
where $\cX$ is a smooth proper algebraic space over 
$B$ and $\varrho_b:\cX_b \ra \cX'_b$
is the minimal resolution of $\cX'_b$ for each $b\in B$.  
However, $\cX$ is constructed by making modifications
of $\cX'$ in formal neighborhoods of the singular fibers,
and hence is not necessarily a scheme.  
Note that blowing up the singularities of $\cX'$
will usually introduce exceptional divisors in the fibers.  
\end{exam}

\begin{defi}
Let $b\in B\setminus S$ be a point of good reduction and 
$
\pi: \cX\ra B
$ 
a good model. An $N$-{\em jet} of 
$\pi$ at $b$ is a section of 
$$
\cX\times_B \Spec(\cO_{B,b}/\mfm^{N+1}_{B,b})\ra 
\Spec(\cO_{B,b}/\mfm^{N+1}_{B,b}).
$$
\end{defi}

Hensel's lemma guarantees that every $N$-jet is 
a restriction of a section of $\hat{\pi}_b$. 
Let $\{b_i\}_{i\in I}$ be a finite set of points of 
good reduction and $j_i$ an $N$-jet of $\pi$ at $b_i$. 
We write $J=\{ j_i\}_{i\in I}$ for the
corresponding collection of $N$-jets.

\subsection{Iterated blowups}
\label{sect:iter}
Let $\pi:\cX\ra B$ be good model of $X$
and $J=\{ j_i\}_{i\in I}$ a finite collection of $N$-jets
at points of good reduction $\{b_i\}$.   
The {\em iterated blowup} associated with $J$
$$\beta(J):\cX(J)\ra \cX$$ 
is obtained by performing 
the following sequence of blowups:
For each $i\in I$ choose a section $\hat{s}_i$ of $\hat{\pi}_{b_i}$ 
with jet $j_i$. 
Now blow up $\cX$ successively $N$ times, where at each stage
the center is the point at which
the proper transform of
$\hat{s}_i$ meets
the fiber over $b_i$.  
Observe that a blowup of $\cX$ centered in the  
fibers of $\pi$ is uniquely determined by the corresponding 
blowup of the completions along those fibers.
Note that at each stage we blow up a smooth point of 
the fiber of the corresponding model and that the 
result does not depend on the order of the $b_i$ or
on the choice of $\hat{s}_i$.

The fiber $\cX(J)_{b_i}$ decomposes into irreducible components
$$
\cX(J)_{b_i}=E_{i,0} \cup  \ldots   \cup E_{i,N}
$$
where 
\begin{itemize}
\item $E_{i,0}$ is the proper transform of $\cX_{b_i}$,
isomorphic to the blowup of $\cX_{b_i}$ at $r_{i,0}:=\hat{s}_i(b_i)$;
\item $E_{i,n}$,  $n=1,\ldots, N-1$, is the blowup of $\bP^d$ 
at $r_{i,n}$, the point where the proper transform 
of $\hat{s}_i$ meets the fiber over $b_i$ of the $n$th blowup;
\item $E_{i,N}\simeq \bP^d$.
\end{itemize}
The intersection $E_{i,n}\cap E_{i,n+1}$ 
is the exceptional divisor $\bP^{d-1}\subset E_{i,n}$
and a proper transform of a hyperplane
in $E_{i,n+1}$, for $n=0,\ldots, N-1$.  
Let $r_i \in E_{i,N}\setminus E_{i,N-1}$ denote
the intersection of $\hat{s}_i$ with $E_{i,N}$.

\begin{figure}[htb]
\centerline{\hskip 7cm\includegraphics[width=\textwidth]{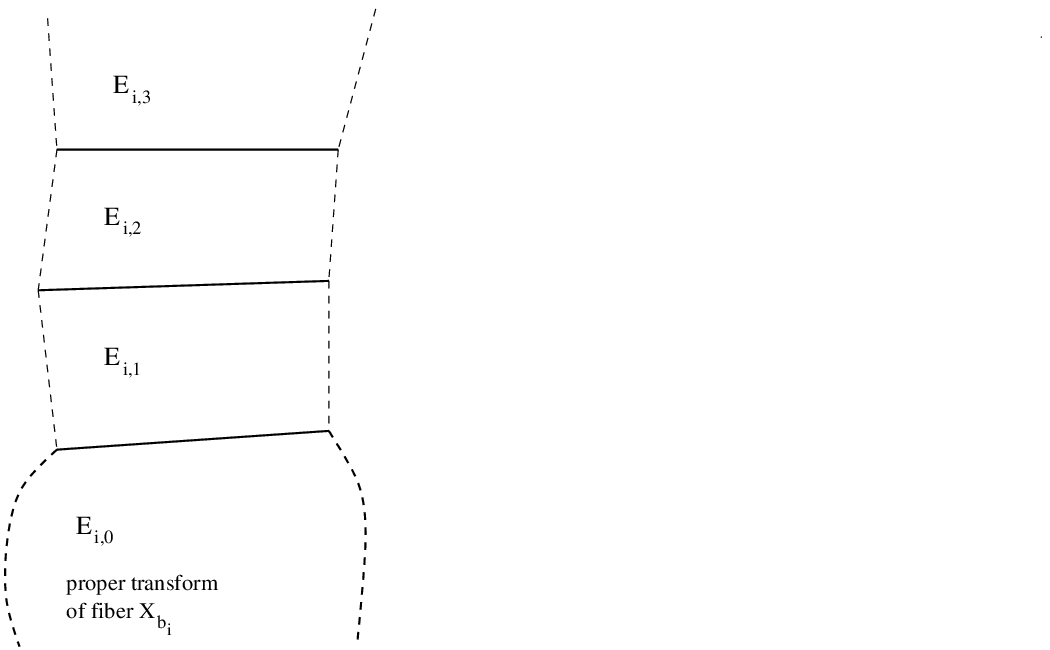}}
\caption{One fiber of the iterated blowup}
\label{fig1}
\end{figure}

For each section $s':\cX(J)\ra B$ the composition $\beta(J)\circ s'$
is a section of $\pi:\cX \ra B$.  Conversely, given a section $s$ of
$\pi:\cX \ra B$, its proper transform $s'$ is the unique section
of $\cX(J)\ra B$ lifting $s$.  
Sections $s'$ of $\cX(J)\ra B$ with $s'(b_i)=r_i$
yield sections of $\cX\ra B$ with
$N$-jet $j_i$ at $b_i$.  
We therefore have the following criterion for weak
approximation in fibers of good reduction:

\begin{prop}
\label{prop:jets}
$X$ satisfies weak approximation in fibers of good reduction
if and only if 
there exists a good model $\pi:\cX \ra B$ with the following
property:  
For each $N\in \bN$ and collection of $N$-jets $J$ 
with corresponding iterated blowup $\cX(J)$, and for
any choice of points $r_i\in E_{i,N}\setminus E_{i,N-1}, i\in I$,
there exists a section $s$ of $\cX(J)\ra B$ with $s(b_i)=r_i$
for each $i\in I$.  
\end{prop}

\section{Rationally connected varieties}
\label{sect:rcv}
We retain the notation introduced in Section \ref{sect:not}.
In particular, the ground field $k$ is algebraically closed 
of characteristic zero. 

\subsection{Terminology and fundamental results}
\label{sect:rat}
Rational connectedness
was introduced in the classification of Fano
varieties \cite{Ca} \cite{KMM}.  However, rationally connected
varieties are now of independent interest:
\begin{defi}[\cite{kollar} IV.3.2]
\label{defi:ko}
A variety $Y$ is {\em rationally chain connected}
(resp. {\em rationally connected}) if there is a 
family of proper and connected curves $g:U\ra Z$
whose geometric fibers have only rational components
(resp. are irreducible rational curves) and a
cycle morphism $g:U\ra Y$ such that 
$$
u^{(2)}:U\times_Z U \ra Y\times Y.
$$
is dominant.
\end{defi}
Our definition of `rationally chain connected' makes sense
for reducible schemes $Y$.

\begin{exam}The class of rationally connected varieties includes
unirational varieties and smooth Fano varieties, \cite{kollar} IV.3.2.6, 
V.2.13.  In particular, smooth hypersurfaces of degree $\le m$
in $\bP^m$ are rationally connected.
\end{exam}

\begin{defi}
Let $Y$ be a smooth algebraic space of dimension $d$ and
$f:\bP^1 \ra Y$ a morphism, so we have an isomorphism
$$
f^*T_Y\simeq \cO_{\bP^1}(a_1) \oplus \ldots \oplus \cO_{\bP^1}(a_d)
$$
for suitable integers $a_1,\ldots,a_d$.
Then $f$ is {\em free} (resp. {\em very free}) if each $a_i\ge 0$
(resp. $a_i\ge 1$).  
\end{defi} 

We recall some key properties:
\begin{itemize}
\item[(1)]Let $Y$ be a proper rationally chain connected variety.  
Then any two closed points are contained in a connected curve
with rational irreducible components, \cite{kollar} IV.3.5.1.
\item[(2)]
If the ground field $k$ is uncountable 
then $Y$ is rationally connected (resp. rationally
chain connected) if any two very general closed points $y_1$
and $y_2$ are contained in an irreducible rational curve
(resp. connected curve with rational irreducible components), 
\cite{kollar} IV.3.6.
\item[(3)] Let $Y$ be a smooth proper rationally connected variety
and $y_1,\ldots,y_m$ points in $Y$.  Then  
there exists a very free morphism $f:\bP^1 \ra Y$
such that $y_1,\ldots,y_m \in f(\bP^1)$.  We may take $f$
to be an immersion if $\dim(Y)=2$ and an embedding if $\dim(Y)\ge 3$,
\cite{kollar} IV.3.9.
\item[(4)]
A smooth variety $Y$ is rationally connected if
it is rationally chain connected, \cite{kollar} IV.3.10.3.
\item[(5)]
Let $\pi:\cY \ra B$ be a proper equidimensional morphism over
an irreducible base.  If the generic fiber of $\pi$ is 
rationally chain
connected then every fiber is rationally {\em chain} connected, 
\cite{kollar} IV.3.5.2.
\item[(6)]
If $\pi:\cY \ra B$ is a smooth morphism then the locus 
$$\{ b\in B: \cY_b \text{ is rationally connected } \}$$
is open, \cite{kollar} IV.3.11.  
\end{itemize}
Since Properties (5) and (6) are local on the base, 
they also hold for good models (which are only assumed to
be algebraic spaces over $B$).

\begin{exam}
\label{exam:cubic}
Property (5) does not guarantee that every
fiber is rationally connected: Consider the
family of cubic surfaces
$$
\cX:=\{(w,x,y,z;t): x^3+y^3+z^3=tw^3 \}\ra \bA^1_t.
$$
The generic fiber is rationally connected but the fiber
$\cX_0$ is a cone over an elliptic curve, which
is not rationally connected.
\end{exam}

\subsection{Producing sections through prescribed points}
Theorem~\ref{theo:great}, when combined with the machinery of
Section~\ref{sect:rat}, has the following
important consequence:

\begin{theo}[\cite{kollar} IV.6.10, \cite{KMM} 2.13]
\label{theo:ghsk}
Let $X$ be a smooth projective
rationally connected variety over the
function field of a curve.  Given a projective model $\pi:\cX \ra B$,
a finite collection of points $\{b_i \}_{i\in I}$ such that each
$\cX_{b_i}$ is smooth, 
and points $x_i \in \cX_{b_i}$, there exists a section
$s:B\ra \cX$ such that $s(b_i)=x_i$, for each $i\in I$.  
\end{theo}

It is natural to wonder whether we can relax the hypothesis that
the fibers $\cX_{b_i}$ be smooth.  For simplicity, assume
that the total space $\cX$ of our model is regular;  this 
can always be achieved by resolving singularities.  Then
for each section $s$, $s(b)\in \cX_b$ is necessarily a smooth point;
otherwise, the intersection multiplicity of the section with
$\cX_b$ would be $>1$.  In light of this, the most optimistic
generalization of Theorem~\ref{theo:ghsk} would be:
\begin{conj}
Let $X$ be a smooth projective
rationally connected variety over the
function field of a curve.  Given a regular model $\pi:\cX \ra B$,
a finite collection of points $\{b_i \}_{i\in I}\subset B$
and smooth points $r_i \in \cX_{b_i}$, there exists a section
$s\,:\,B\ra \cX$ such that $s(b_i)=r_i$ for each $i\in I$.  
\end{conj}

Applying this to the iterated blowups as described 
in Proposition~\ref{prop:jets},
we obtain: 
\begin{conj}
\label{conj:main}
A smooth rationally connected variety
over the function field of a curve 
satisfies weak approximation.  
\end{conj}

\

We outline the main issues in the proof of Theorem~\ref{theo:main};
details are given in Section~\ref{sect:main-proof}. 
By Proposition~\ref{prop:jets}, we are reduced to proving the
existence of a section passing through specific smooth points $r_i$
of {\em singular} fibers of the iterated blow-up.
Theorem~\ref{theo:ghsk} does not immediately imply this, but it does
guarantee a section $\sigma$ passing through some point 
$x_i$ of each of these fibers.  
Property (5) of rationally connected varieties from
Section~\ref{sect:rat} guarantees the existence
of some chain $T_i$ of rational curves in the corresponding fiber
joining $x_i$ and $r_i$. 
The difficulty is to choose these so that
$C:=\sigma(B) \cup_{i\in I} T_i$ deforms to a 
section containing $r_i$, for each $i\in I$. 
In particular, it is necessary that $C$ intersect the 
fibral exceptional divisor $E_{i,N}$ containing $r_i$ in one point
and not intersect the other components of the corresponding fiber;
this constrains the homology class of $T_i$.  
Furthermore, we must describe each $T_i$ explicitly so the deformation
space of $C$ can be analyzed.

\section{Deformation theory}
We continue to work over an algebraically closed ground field
of characteristic zero.
In this section, a {\em curve} is a connected reduced scheme
of dimension one.  
\subsection{Combs}
Recall the {\em dual graph} associated with a nodal curve $C$:
Its vertices are indexed by the irreducible components of $C$ 
and its edges are indexed by 
the intersections of these components. 
\begin{defi}
A projective nodal curve $C$ is {\em tree-like} if
\begin{itemize}
\item{each irreducible component of $C$ is smooth;}
\item{the dual graph of $C$ is a tree.}
\end{itemize}
\end{defi}

We shall require a slight generalization of the standard notion
of a comb (cf. \cite{kollar}):
\begin{defi}
A {\em comb with $m$ broken teeth} is a projective nodal curve $C$ with
$m+1$ subcurves
$D, T_1, \ldots, T_m$ such that 
\begin{itemize}
\item{$D$ is smooth and irreducible;}
\item{$T_{\ell}\cap T_{\ell'}=\emptyset$, for all $\ell\neq \ell'$;}
\item{each $T_{\ell}$
meets $D$ transversally
in a single point; and }
\item{each $T_{\ell}$ is a chain of $\bP^1$'s.}
\end{itemize}
Here $D$ is called the {\em handle} and the
$T_{\ell}$ the {\em broken teeth}.
\end{defi}

\begin{figure}[h]
\centerline{\includegraphics{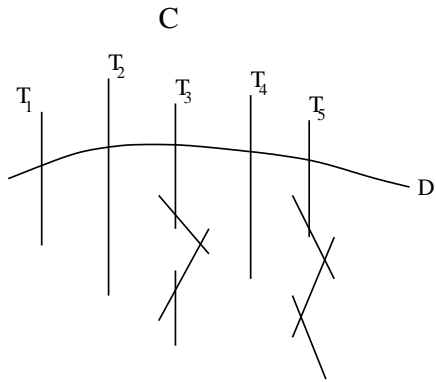}}
\caption{A comb with five broken teeth}\label{fig3}
\end{figure}

\subsection{Vector bundle lemmas}
\label{sect:vblem}
Let $C$ be a smooth curve and $\cV$ a vector
bundle on $C$.  Given a collection of distinct
points $\mfq=\{q_1,\ldots,q_m\}\subset C$ and one-dimensional
subspaces of the fibers $\xi_{q_{\ell}} \subset \cV_{q_{\ell}},\ell=1,\ldots,m$,
there exists a rank-one subbundle $\cL \subset \cV$ with fiber at
$q_{\ell}$ equal to $\xi_{q_{\ell}}$, $\ell=1,\ldots,m$.   
The extension
$$
0 \longrightarrow \cL \longrightarrow \cL\otimes  \cO_C(\mfq)
\longrightarrow \cL \otimes (\cO_C(\mfq))_{\mfq} \longrightarrow 0
$$
induces
$$
0 \longrightarrow \cV \longrightarrow \cV' \longrightarrow Q 
\longrightarrow 0,
$$
where $Q:=(\cL\otimes \cO_C(\mfq))_{\mfq}$ 
is supported on $\mfq$ and has
length one at each $q_{\ell}$.  
This extension depends on the $q_{\ell}$ and $\xi_{q_{\ell}}$ but not 
on $\cL$.  The saturation of $\cL$ in $\cV'$ is isomorphic
to $\cL\otimes \cO_C(\mfq)$. 
\begin{lemm} \label{lemm:sat}
Retain the notation introduced above.  
A subbundle $\cU \subset \cV$ is also a subbundle of $\cV'$ if 
$\xi_{q_{\ell}}\cap \cU_{q_{\ell}}=0$ for each $\ell$.  
\end{lemm}
\begin{proof}
We have exact sequences
\begin{eqnarray*}
& & 0 \ra \Hom(\cL\otimes \cO_C(\mfq),\cV) \ra \Hom(\cL,\cV) \ra
\Ext^1(Q,\cV) \\
& & 0 \ra \iHom(\cL\otimes \cO_C(\mfq),\cV) \ra \iHom(\cL,\cV) \ra
\iExt^1(Q,\cV)  \ra 0
\end{eqnarray*}
and the extension class $\eta_{\cV'} \in \Ext^1(Q,\cV)$
is the image of the inclusion $\cL \hookrightarrow \cV$ under 
the connecting homomorphism.  Since $\cL$ is saturated in
$\cV$ at $q_{\ell}$, $\eta_{V'}$ localizes to a nonzero
element of $\iExt^1(Q,\cV)_{q_{\ell}}$
for each $\ell$.  

Since $\cU\subset \cV$ is a subbundle, $\cV/\cU$ has no torsion and
thus is locally free.  
The class $\eta_{\cV'}$ naturally induces an
extension
$$0 \longrightarrow \cV/\cU \longrightarrow \cV'/\cU \longrightarrow Q \ra 0,$$
classified by $\eta_{\cV'/\cU}\in \Ext^1(Q,\cV/\cU)$,
the image of the
composition 
$$
\cL \hookrightarrow \cV \twoheadrightarrow \cV/\cU
$$
under the connecting homomorphism.  Our hypothesis guarantees
that $\cL$ is a subbundle of $\cV/\cU$ near $q_{\ell}$, hence
$\eta_{\cV'/\cU}$ localizes to a nonzero element of 
$\iExt^1(Q,\cV/\cU)_{q_{\ell}}$ for each $\ell$.
It follows that $\cV'/\cU$ is torsion-free, so $\cU$ is a subbundle.
\end{proof}

\begin{lemm}[\cite{GHS}, Lemma 2.5]
\label{lemm:vb}
Retain the notation introduced above.  Fix an integer $N$ and 
a vector bundle $\cV$.  Then there exist points 
$q_1, \ldots, q_m$ and one-dimensional subspaces $\xi_{q_\ell}\subset \cV_{q_{\ell}}$, $\ell=1,\ldots,m$, such that
$$H^1(\cV'\otimes \cO_C(-w_1-\ldots-w_N))=0$$
for any points $w_1,\ldots,w_N \in C$.  
\end{lemm}

Our next lemma is well known (cf. \cite{NS}, Section 2) but we
provide a proof for the convenience of the reader:
\begin{lemm}
\label{lemm:tr-like}
Let $C$ be a tree-like curve and $\cV$ a vector bundle on $C$. 
If for each irreducible component $C_{\ell}$ of $C$
the restriction $\cV\otimes \cO_{C_{\ell}}$ is globally generated
then $\cV$ is globally generated.
Furthermore, 
$$H^1(C,\cV)\ra \oplus_{\ell} H^1(\cV\otimes \cO_{C_{\ell}})$$ 
is an isomorphism.  
\end{lemm}

\begin{proof}We do induction on the number of
irreducible components;  the case of one component
is trivial.  Otherwise, express $C$ as a union $D\cup D^c$,
where $D$ is irreducible with connected complement in $C$
and $D^c=\overline{C\setminus D}$ is tree-like.  Let $q$
be the node of $C$ joining $D$ and $D^c$,
$$
g:C':=D\amalg D^c\ra C
$$
the partial normalization of $C$ at $q$, and  
$r,r^c$ the points of $C'$
with $g(r)=g(r^c)=q$.    
The descent data for $\cV$ consist of the pullback $g^*\cV$
and an isomorphism
$$\phi:(g^*\cV)_r \ra (g^*\cV)_{r^c}$$
induced by identifications of fibers 
$$
(g^*\cV)_r\simeq \cV_q \simeq (g^*\cV)_{r^c}.
$$
Recall the exact sequences relating the cohomology 
of $\cV$ and $g^*\cV$:
$$
0 \ra  H^0(C,\cV) \ra
H^0(C',g^*\cV) \ra
(g^*\cV)_r\oplus (g^*\cV)_{r^c}
  \stackrel{(-\phi,\mathrm{Id})}{\ra}
  (g^*\cV)_{r^c} \ra 0 
$$
$$
H^0(C',g^*\cV)  \ra
(g^*\cV)_r\oplus (g^*\cV)_{r^c}  \ra
H^1(C,\cV) \ra H^1(C',g^*\cV) \longrightarrow 0.
$$

By the inductive hypothesis, $g^*\cV$ is globally generated 
on $C'$.  Since $r$ and $r^c$ are on different connected
components of $C'$, the second exact sequence guarantees
that $H^1(C,\cV)\ra H^1(C',g^*\cV)$ is injective;
the cohomology statement follows.
Since $g^*\cV$ is globally generated, for each section over
$D$ there exists a section over $D^c$ compatible under
the isomorphism $\phi$, and vice versa.  These compatible
pairs of sections descend to elements of $H^0(C,\cV)$.
Thus given $p\in D$ and $v\in \cV_p$, a section 
$t\in H^0(D,\cV\otimes \cO_D)$ with $t(p)=v$ 
extends to a section over $C$.  
\end{proof}

\subsection{Analysis of normal bundles}
\label{sect:nb}
We describe
the normal bundle of a nodal curve immersed in
a smooth algebraic space.  Our main references are
Section 2 of \cite{GHS} and Section 6 of \cite{AK}.
See \cite{Art69} and \cite{kollar} I.5
for foundational results on Hilbert `schemes' of algebraic spaces,
\cite{kollar} II.1
for applications to morphisms of curves into spaces,
and \cite{dJS} for an extension of
Theorem~\ref{theo:great} to 
positive characteristic using this machinery.

If $C$ is a nodal curve imbedded into a smooth space $Y$
then $\cN_{C/Y}$ is defined as the dual to the kernel of the restriction
homomorphism of K\"ahler differentials
$\Omega^1_Y\otimes \cO_C \twoheadrightarrow \Omega^1_C;$
a local computation shows this is locally free.
First order deformations of $C\subset Y$ are given by
$H^0(C,\cN_{C/Y})$;  obstructions are given by $H^1(C,\cN_{C/Y})$.  

Let $D$ be the union of irreducible components of
$C$ and $\mfq=\{q_1,\ldots,q_m\}$ the locus
where $D^c:=\overline{C\setminus D}$ meets $D$.
At a node of $C$, the tangent cone is a union of two one-dimensional
subspaces, the tangents to the transverse branches.
The tangent to $D^c$ at $q_{\ell}$ yields a one-dimensional subspace
$$\xi_{q_{\ell}}\subset (\cN_{D/Y})_{q_{\ell}}.$$
As in Section \ref{sect:vblem}, these induce a natural extension
$\cN'_{D/Y}$ of $\cN_{D/Y}$, which coincides with
the restriction to $D$ of the normal bundle to $C$ in $Y$
$$
0 \lra \cN_{D/Y} \lra \cN_{C/Y}\otimes \cO_D \lra Q \lra 0.
$$
Here $Q$ is a torsion sheaf supported on $\mfq$, with length one
at each point $q_{\ell},\ell=1,\ldots,m$.  Sections of 
$\cN_{C/Y}\otimes \cO_D$ can be interpreted as sections of $\cN_{D/Y}$
with poles at the $q_{\ell}$ in the directions $T_{q_{\ell}}D^c$.  

We shall need a slight generalization:  We continue to assume that
$C$ is a nodal curve and $Y$ is nonsingular.  Let $f:C \ra Y$ denote
a closed immersion whose image is a nodal curve.  The 
restriction homomorphism
$f^*\Omega^1_Y \rightarrow \Omega^1_C$
is surjective and the dual to its kernel is still locally free.
This is denoted $\cN_f$ and coincides with $\cN_{C/Y}$ when $f$
is an embedding.  
First order deformations of $f:C\ra Y$ are given by
$H^0(C,\cN_f)$;  obstructions are given by $H^1(C,\cN_f)$.  
This set-up differs from the standard deformation theory of morphisms
in that we ignore reparametrizations of $C$.  
When $D$ is a union of irreducible
components of $C$ as above, then the analogous extension
takes the form (cf. Lemma 2.6 of \cite{GHS}):
$$
0 \lra \cN_{f|_D} \lra \cN_f \otimes \cO_D \lra Q \lra 0.
$$

This analysis gives the following infinitesimal smoothing criterion:
(cf. Lemma 2.6 of \cite{GHS}):
\begin{lemm} \label{lemm:smooth}
Retain the above notation.  A first
order deformation $t\in H^0(C,\cN_f)$ smooths the node
$q_{\ell}$ if the restriction
$t|_D \in \cN_f\otimes \cO_D$ is nonzero in $Q_{q_{\ell}}$.
\end{lemm}

\begin{prop}
\label{prop:diag}
Let $Y$ be a smooth algebraic space, $E\subset Y$ a smooth subspace of
codimension one, $C$ a nodal curve, $D\subset C$ a union of irreducible
components of $C$, $D^c=\overline{D\setminus C}$, and
$\mfq=\{q_{\ell}\}=D\cap D^c$.  
Let $f:C \rightarrow Y$ be an immersion with image a nodal
curve so that $f(D)\subset E$, $f(D\setminus C)\subset Y\setminus E$,
and $f(D^c)$ is transverse to
$E$ at each point of $f(\mfq)$.
If $g:D\ra E$ is the resulting immersion of $D$
into $E$ then
$\cN_g$ is saturated in $\cN_f\otimes \cO_D$
and we have the following diagram:

\centerline{
\xymatrix{
        &                         &  0 \ar[d]               & 0\ar[d]               &  \\
0\ar[r] &\cN_g\ar[r]\ar@{=}[d]& \cN_{f|_D}\ar[r]\ar[d]   & g^*\cN_{E/Y} \ar[d]\ar[r] & 0 \\
0\ar[r] &\cN_g\ar[r]          & \cN_f\otimes \cO_D\ar[r]\ar[d]& g^*\cN_{E/Y}\otimes \cO_D(\mfq)\ar[d]\ar[r] & 0 \\
        & 0\ar[r]                         &  Q\ar[r]\ar[d]      & g^*\cN_{E/Y}\otimes (\cO_D(\mfq))_{\mfq}\ar[d] \ar[r]& 0\\
        &                         &   0                     &   0                                        & 
}}
   
\end{prop}

\begin{proof}
For each $\ell$ the composition
$$T_{q_{\ell}}D^c \ra T_{q_{\ell}}Y \ra 
T_{q_{\ell}}Y/ T_{q_{\ell}}D$$
determines a one-dimensional subspace
$\xi_{\ell}\subset (\cN_{f|_D})_{q_{\ell}}$.
The transversality hypothesis implies
 $(\cN_g)_{q_{\ell}}\cap \xi_{q_{\ell}}=0$ in $(\cN_{f|_D})_{q_{\ell}}$. 
Using Lemma \ref{lemm:sat}, we conclude that $\cN_{f|_D}$ is
a subbundle in $\cN_f\otimes \cO_D$, so the quotient 
$$R=\left(\cN_f\otimes \cO_D\right) / \cN_{f|_D}$$
is locally free.  This sheaf arises as an extension
\begin{equation}
\label{eq:extend}
0 \longrightarrow g^*\cN_{E/Y} \longrightarrow R \longrightarrow 
Q \longrightarrow 0,
\end{equation}
where $Q$ is supported on $\mfq$ with length
one at each $q_{\ell}$.  We may therefore identify extension
(\ref{eq:extend}) with  the tensor product of $g^*\cN_{E/Y}$ with
$$
0 \longrightarrow \cO_D \longrightarrow \cO_D(\mfq) 
\longrightarrow \left(\cO_D(\mfq)\right)_{\mfq} \longrightarrow 0.
$$
\end{proof}

\begin{prop}
\label{prop:tree}
Let $C$ be a tree-like curve, $Y$ a smooth algebraic space, and
$f:C \ra Y$ an immersion with nodal image.
Suppose that for each irreducible component $C_\ell$ of $C$, 
$H^1(C_{\ell},\cN_f\otimes \cO_{C_{\ell}})=0$ and 
$\cN_f\otimes \cO_{C_{\ell}}$
is globally generated.  
Then $f:C \ra Y$ deforms to an immersion of a smooth curve
into $Y$.
   
Suppose furthermore that
$\mfp=\{p_1,\ldots,p_w \}\subset C$ 
is a collection of smooth points 
such that for each component $C_{\ell}$
$H^1(\cN_f\otimes \cO_{C_{\ell}}(-\mfp))=0$ 
and the sheaf $\cN_f\otimes \cO_{C_{\ell}}(-\mfp)$  
is globally generated.
Then $f:C \ra Y$ deforms to an immersion of a smooth curve
into $Y$ containing $f(\mfp)$.
\end{prop} 

\begin{proof}
Our argument is similar to the constructions of Section 2 of
\cite{GHS} and Lemma 65 of \cite{AK}.  
Lemma~\ref{lemm:tr-like} implies that $H^1(C,\cN_f)=0$ and $\cN_f$
is globally generated.  Hence 
the space of maps is unobstructed and every first-order
deformation of $f$ lifts to an actual deformation.  
Global generation implies the existence of 
$t\in H^0(C,\cN_f)$ so that,
for each component $C_\ell$, the image of $t$ in
$$
H^0(C_{\ell},Q(C_{\ell})), \quad 
Q(C_{\ell}):=\left(\cN_f\otimes \cO_{C_{\ell}}\right)/\cN_{f|_{C_{\ell}}}
$$
is nonzero at each point of the support of $Q(C_{\ell})$.
The first-order deformation $t$ smooths each node of $C$ by 
Lemma~\ref{lemm:smooth}.

For the second part, consider those maps with image containing
$f(\mfp)$.  Our cohomology assumption 
guarantees that this space is unobstructed;
in addition, $\cN_f\otimes \cO_C(-\mfp)$ is globally generated.  
Hence there exists  a $u\in H^0(\cN_f\otimes \cO_C(-\mfp))$ so that,
for every component $C_\ell$, the image of $u$ in $Q(C_{\ell})$
is nonzero at each point of its support.
Note that 
$$
(\cN_f\otimes \cO_{C_{\ell}}(-\mfp))/
(\cN_{f|_{C_{\ell}}}\otimes \cO_{C_{\ell}}(-\mfp))=
(\cN_f\otimes \cO_{C_{\ell}})/\cN_{f|_{C_{\ell}}},
$$
since the quotient is a torsion sheaf with support disjoint from $\mfp$. 
Hence the first-order deformation $u$ smooths each node of $C$ and 
contains $f(\mfp)$.
\end{proof}

\section{Proof of the main theorem}
\label{sect:main-proof}

The theorem is well known when $d=\dim(X)=1$.  The
only smooth proper rationally connected curve is $\bP^1$,
which satisfies weak approximation.  We may therefore assume
$d\ge 2$.

Recall the set-up of Proposition~\ref{prop:jets}:
It suffices to show that for
each integer $N$, finite set $\{b_i\}_{i\in I}\subset B\setminus S$,
and collection of $N$-jet data $J$ supported in the
fibers over $\{b_i\}_{i\in I}$, there exists
a section in the iterated blowup $\cX(J)$ passing through
prescribed points $r_i \in E_{i,N}\setminus E_{i,N-1}$.  

We proceed by induction on $N$;  the base case $N=0$ is essentially
Theorem~\ref{theo:ghsk}.  However, our assumptions are slightly weaker:
We are not assuming $\cX$ is a scheme.  

The total space of $\cX$ is smooth along the fibers $\cX_{b_i}, i\in I$,
so we may resolve the singularities of $\cX$ without altering these fibers.
Theorem~\ref{theo:great} gives a section $\sigma$ of $\pi$.
Let $q_i=\sigma(b_i)$ for each $i\in I$;  let $I'\subset I$
(resp. $I''\subset I$) denote those indices with $q_i\neq r_i$ 
(resp. $q_i =r_i$).  

We shall construct a comb $C$ with handle $\sigma(B)$ and smooth
teeth $T_1,\ldots,T_m$ and an immersion $f:C \ra \cX$
with nodal image so that:
\begin{itemize}
\item{the $T_{\ell}$ are free rational curves
in distinct smooth fibers of $\pi$;}
\item{for each $i\in I'$, there is a tooth $T_i$
containing $r_i$ as a smooth point;}
\item{let $\mfr$ denote the sum of the points of $C$
mapping to the $r_i,i\in I$;  then the restriction of 
$\cN_f\otimes \cO_C(-\mfr)$
to each irreducible component of $C$ 
is globally generated and has no higher cohomology.}
\end{itemize}
We emphasize that $f$ can be taken to be an
embedding if $d>2$.  

Proposition~\ref{prop:tree} implies $f:C\ra \cX$
admits a deformation $\tilde{f}:\tilde{C}\ra \cX$, 
where $\tilde{C}$ is smooth and $\tilde{f}(\tilde{C})$ contains 
each of the $r_i$.  Since all the teeth $T_1,\ldots,T_m$
are contained in fibers of $\cX\ra B$, 
$C$ intersects the generic fiber
in one point.  Thus the deformed curve $\tilde{C}$ also
meets the generic fiber in one point and hence is a section
of $\cX \ra B$.  

Here are the details of the construction.  For each $\ell$
$$f_{\ell}:T_{\ell} \ra \cX_{b_{\ell}}$$
is a free rational curve with nodal image,
so that $\sigma(b_{\ell})\in f_{\ell}(T_{\ell})$ as a smooth
point.  For $\ell=1,\ldots,|I'|$ we choose these so that
$r_i$ is contained in the image as a smooth point.  
For $\ell=|I'|+1,\ldots,m$, we choose these in generic fibers
of good reduction with generic tangent directions
$\xi_{\ell} \subset T_{\sigma(b_{\ell})}\cX_{b_{\ell}}$
satisying the hypotheses of Lemma~\ref{lemm:vb}, 
so that the extension
$$
0 \lra \cN_{\sigma(B)} \lra \cN_f\otimes\cO_{\sigma(B)} \lra Q(\sigma(B)) \lra 0
$$
is globally-generated and has no higher cohomology, even
after twisting by $\cO_{\sigma(B)}(-\sum_{i\in I''}r_i)$.

\

We next address the inductive step.  
Let $J'$ denote order-$(N-1)$ truncation
of $J$, i.e., if 
$$
j_i: \Spec(\cO_{B,b_i}/\mfm_{B,b_i}^{N+1}) \ra 
\cX\times_B\Spec(\cO_{B,b_i}/\mfm_{B,b_i}^{N+1})
$$
then
$$j'_i=j_i|\Spec(\cO_{B,b_i}/\mfm_{B,b_i}^N).$$
The inductive hypothesis applied to $J'$ guarantees the
existence of a section $s':B\ra \cX(J')$ passing through arbitrary
points 
$$
r'_i\in (E'_{i,N-1}\setminus E'_{i,N-2})\subset \cX(J').
$$  
Specifically,
we choose $s'$ so that it has jet data $J'$ over the
points $\{b_i\}_{i\in I}$.  Let $\sigma:B \ra \cX(J)$ denote
the proper transform of $s'$ in $\cX(J)$.  By construction,
$\sigma$ meets $\cX(J)_{b_i}$ in a point 
$q_i\in E_{i,N}\setminus E_{i,N-1}$ for
each $i$.  Our goal is to find a section $s:B \ra \cX(J)$ 
such that for each $i\in I$ and $r_i\in E_{i,N}\setminus E_{i,N-1}$,
$s(b_i)=r_i$.  
Again,  $I'\subset I$
(resp. $I''\subset I$) denotes those indices with $q_i\neq r_i$ 
(resp. $q_i =r_i$).  

Next, we construct a comb $C$ with handle $\sigma(B)$ and 
broken teeth $T_1,\ldots,T_m$ and an immersion $f:C \ra \cX(J)$
with nodal image so that:
\begin{itemize}
\item{for each $i\in I'$,
there is a broken tooth $T_i$
mapped to $\cX(J)_{b_i}$ and containing $r_i$;}
\item{$C$ is smoothly embedded at $r_i$ for each $i\in I'$, so there is a unique
component $T_{i,N}\subset C$ containing $r_i$;}
\item{the remaining broken teeth $T_{|I'|+1},\ldots,T_m$
are free rational curves contained
in generic fibers of $\cX(J)\ra B$ of good reduction;} 
\item{the restriction of $\cN_f\otimes \cO_C(-\mfr)$ 
to each irreducible component
is globally generated and has no higher cohomology.}
\end{itemize}
Again, $f$ can be taken to be an
embedding if $d>2$.  

Proposition~\ref{prop:tree} implies that $f:C\ra \cX(J)$
admits a deformation $\tilde{f}:\tilde{C}\ra \cX(J)$, 
where $\tilde{C}$ is smooth and $\tilde{f}(\tilde{C})$ contains 
each of the $r_i$.  The image $\tilde{f}(\tilde{C})$
is the desired section of $\cX(J) \ra B$.  

We start by describing the teeth $T_i$ with $i\in I'$.
Recall from Section \ref{sect:iter}
that $E_{i,N}\simeq \bP^d$ and $E_{i,N}\cap E_{i,N-1}$ is a hyperplane
section in this $\bP^d$.  Let $T_{i,N}$ denote the unique line joining
$r_i$ to $q_i=\sigma(b_i)$;  let $q_{i,N-1}$ denote the intersection of
this line with $E_{i,N-1}$.   We have
$$E_{i,N-1}\simeq \mathrm{Bl}_{r_{i,N-1}}\bP^d$$
with exceptional divisor $E_{i,N-1} \cap E_{i,N}\simeq \bP^{d-1}$;
there is a unique line in $\bP^d$ containing $r_{i,N-1}$ whose
proper transform $T_{i,N-1}\subset E_{i,N-1}$ meets $q_{i,N-1}$.
Let $q_{i,N-2}$ denote the intersection of this line with $E_{i,N-2}$.  
Continuing in this way, we obtain a sequence of embedded 
smooth rational curves 
$$T_{i,n}\subset E_{i,n}\simeq \mathrm{Bl}_{r_{i,n}}\bP^d, \quad n>0,$$
each the proper transform of a line meeting $r_{i,n}$.
Let $q_{i,0}$ denote the intersection of $T_{i,1}$ with  
$E_{i,0}$, which is a point in the exceptional divisor of
$$E_{i,0}=\mathrm{Bl}_{r_{i,0}}\ra \cX_{b_i}.$$
Let $g_{i,0}:T_{i,0} \rightarrow E_{i,0}$ be a free 
rational curve, immersed so that the image is a nodal curve, 
with $q_{i,0}\in g_{i,0}(T_{i,0})$
as a smooth point.  Property (3) of rationally connected varieties
gives such a curve;  $g_{i,0}$ can be taken to be an embedding when $d>2$.  
Let $f_{i,0}$ denote the composition of $g_{i,0}$ with the
inclusion $E_{i,0}\subset \cX(J)$, and 
$$f_i:T_i=T_{i,0}\cup \ldots \cup T_{i,N}\lra \cX(J)$$
the resulting map of the broken tooth into $\cX(J)$.

\begin{figure}[h]
\centerline{\hskip7cm \includegraphics{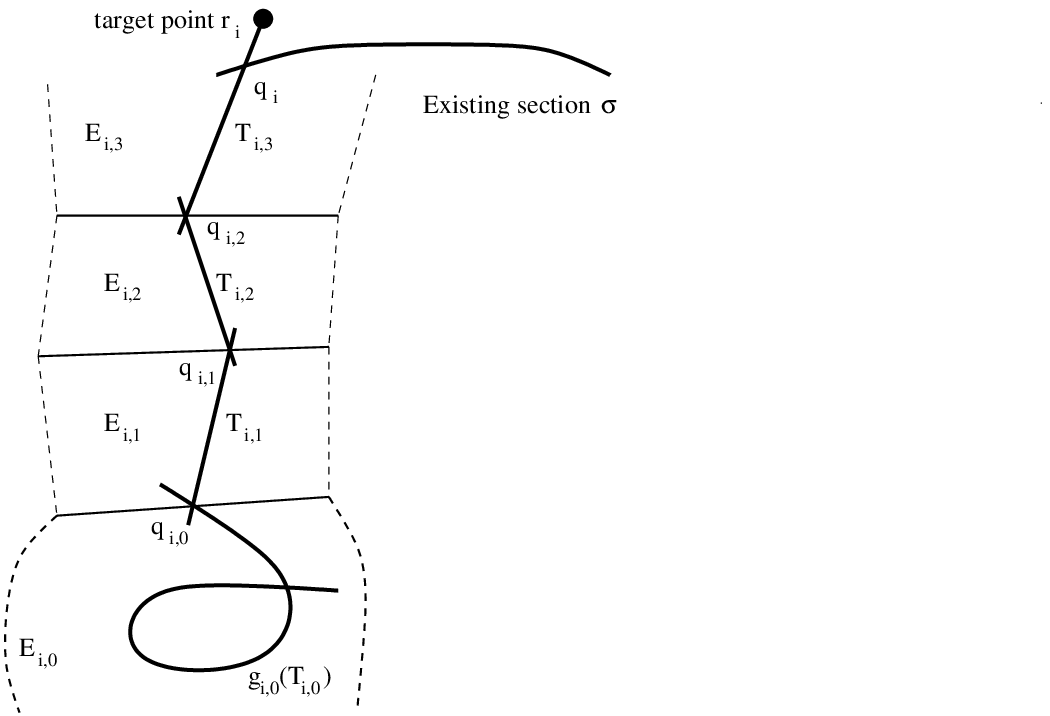}}
\caption{Attaching broken teeth and moving the section}\label{fig2}
\end{figure}

We verify the normal bundle conditions for the components
of these teeth.  First, observe that
\begin{eqnarray*}
\cN_{T_{i,n}/E_{i,n}}&\simeq&\begin{cases}
  \cO_{\bP^1}(+1)^{d-1} & \text{ for } n=N \\
  \cO_{\bP^1}^{d-1} & \text{ for }     n=1,\ldots,  N-1\\
  \end{cases} \\
\cN_{g_{i,0}} &\simeq &
  \oplus_{w=1}^{d-1} \cO_{\bP^1}(a_w), \quad a_w \ge 0;
\end{eqnarray*}
the $a_w$ are nonnegative because $\cN_{g_{i,0}}$ is a quotient
of $g_{i,0}^*T_{E_{i,0}}$, which is nonnegative.  
Fibers of $\cX(J)\ra B$ restrict to the zero divisor on each $T_{i,n}$
and $\sum_{n=0}^N E_{i,n}$ is equivalent to the class of a fiber,
hence
$$
E_{i,n}|T_{i,n}=(-\sum_{n'\neq n} E_{i,n'})|T_{i,n}.
$$
It follows that 
$$
\cN_{E_{i,n}/\cX(J)}\otimes \cO_{T_{i,n}} =
\cO_{\cX(J)}(E_{i,n})\otimes \cO_{T_{i,n}} \simeq
\begin{cases}
  \cO_{\bP^1}(-1) & \text{ for } n=N \\
  \cO_{\bP^1}(-2) & \text{ for } n=1, \ldots, N-1\\
  \cO_{\bP^1}(-1) & \text{ for } n=0.
  \end{cases} 
$$
For $n>0$ we have the exact sequence
$$
0 \lra \cN_{T_{i,n}/E_{i,n}} \lra \cN_{T_{i,n}/\cX(J)} \lra 
\cN_{E_{i,n}/\cX(J)}\otimes \cO_{T_{i,n}} \lra 0,
$$
which splits in our situation.  Therefore, we find
$$
\cN_{T_{i,n}/\cX(J)}\simeq 
\begin{cases}
\cO_{\bP^1}(+1)^{d-1}\oplus \cO_{\bP^1}(-1) & \text{ for } n=N \\
  \cO_{\bP^1}^{d-1}\oplus \cO_{\bP^1}(-2) & \text{ for } n=1,\ldots, N-1
  \end{cases}.
$$
For $n=0$ we have
$$
0 \lra \cN_{g_{i,0}} \lra \cN_{f_{i,0}} \lra 
g_{i,0}^*\cN_{E_{i,0}/\cX(J)} \lra 0,
$$
which implies
$$\cN_{f_{i,0}}\simeq 
  \oplus_{w=1}^{d-1} \cO_{\bP^1}(a_w)\oplus \cO_{\bP^1}(-1).$$

On first examination, the negative summands would make it hard to satisfy
the hypotheses of Proposition~\ref{prop:tree}.  However, the nodes
in the broken teeth give enough positivity to overcome the negative
factors.  We use Proposition~\ref{prop:diag} to analyze
the relationship between the normal bundles to the $T_{i,n}$ and
the restriction to the normal bundle of the broken comb to these
components.  When $n>0$, we have an exact sequence
$$
0 \lra \cN_{T_{i,n}/\cX(J)} \lra \cN_f\otimes \cO_{T_{i,n}}
\lra Q(T_{i,n}) \lra 0,
$$
where $Q(T_{i,n})$ is a torsion sheaf, supported at the nodes of 
$C$ on $T_{i,n}$.  
However, the positive summands of $\cN_{T_{i,n}/\cX(J)}$ are saturated
in $\cN_f$;  only the negative summand fails to be saturated.
When $n=1,\ldots,N-1$, the negative summand is isomorphic 
to $\cO_{\bP^1}(-2)$ and
$Q(T_{i,n})$ has length two, so the saturation 
is $\cO_{\bP^1}$.  When $n=N$, the negative summand is isomorphic to
$\cO_{\bP^1}(-1)$ and $Q(T_{i,N})$ has length two and
support $\{r_i,q_{i,N-1}\}$, so the saturation
is $\cO_{\bP^1}$.  
When $n=0$ we have
$$
0 \lra \cN_{f_{i,0}} \lra \cN_{f_i}\otimes \cO_{T_{i,0}}
\lra Q(T_{i,0}) \lra 0,
$$
with $Q(T_{i,0})$ of length one and supported at $q_{i,0}$.  
The negative summand of $\cN_{f_{i,0}}$ is isomorphic to 
$\cO_{\bP^1}(-1)$, so the extension above induces
$$
0 \ra \cO_{\bP^1}(-1) \ra \cO_{\bP^1} \ra Q(T_{i,0}) \ra 0,
$$
i.e., the saturation of the negative factor is $\cO_{\bP^1}$.  

To summarize, we have shown
$$
\cN_f\otimes \cO_{T_{i,n}} =\begin{cases}
  \cO_{\bP^1}(+1)^d & \text{ for } n=N \\
  \cO_{\bP^1}^d & \text{ for } n=1,\ldots, N-1 \\
 \oplus_{w=1}^{d-1} \cO_{\bP^1}(a_w)\oplus \cO_{\bP^1}, \quad a_w\ge 0&
          \text{ for } n=0
  \end{cases},
$$   
so the hypotheses of  Proposition~\ref{prop:tree}
hold for the broken teeth $T_i$.

For $\ell=|I'|+1,\ldots,m$ we take 
$$f_{\ell}:T_{\ell} \ra \cX_{b_{\ell}}$$
to be free rational curves, immersed in generic fibers of
good reduction so that the images are nodal,
with $\sigma(b_{\ell})\in f_{\ell}(T_{\ell})$ as a smooth
point.  We choose these with generic tangent directions
$\xi_{\ell} \subset T_{\sigma(b_{\ell})}\cX_{b_{\ell}}$
so that Lemma~\ref{lemm:vb} guarantees 
$\cN_f\otimes \cO_{\sigma(B)}(-\sum_{i\in I''}r_i)$
is globally-generated and has 
no higher cohomology.

\bibliographystyle{smfplain}
\bibliography{wa3}

\end{document}